\input amstex
\documentstyle{vis_n}
\magnification 1000
\NoBlackBoxes
\TagsOnRight
\tolerance2500
\footline={\hfill}
\vsize20.5 true cm
\hsize13.0 true cm

\bigskip
\topmatter
\title
ON H-CLOSED PARATOPOLOGICAL GROUPS
\endtitle
\author
Oleksandr RAVSKY
\endauthor
\thanks
\noindent\kern-12pt \copyright \ \ { Ravsky Oleksandr}, 2002
\endthanks
\abstract
A Hausdorff paratopological group is it H-closed if it is closed in every
Hausdorff paratopological group containing it as a paratopological subgroup.
Obtained a criterion when abelian topological group is H-closed
and for some classes of abelian paratopological groups are obtained
simple criteria of H-closedness.

\endabstract
\date {\rm 
Ivan Franko National University of Lviv, 
1 Universitetska  str. 79000 Lviv,  Ukraine 
}
\enddate

\endtopmatter

\abovedisplayskip 6pt plus 3pt minus 1pt
\belowdisplayskip 6pt plus 3pt minus 1pt

\vskip-20pt

\def\inte{\operatorname{int}}
\def\ol{\overline}
\def\0{\varnothing}
\def\bs{\backslash}

\def\N{\Bbb N}
\def\Z{\Bbb Z}
\def\Q{\Bbb Q}

All topological spaces considered in this paper are Hausdorff,
if the opposite is not stated. We shall use the following notations. Let
$A$ be a subset of a group and $n$ be an integer. Put
$A^n=\{a_1a_2\cdots a_n:a_i\in A\}$ and $nA=\{a^n:a\in A\}$. For a
group topology $\tau$ the closure of set $A$ we define as $\ol A^{\tau}$ and
the base of the unit as $\Cal B_\tau$.

A topological space $X$ of a class $C$ of topological spaces is
{\it C-closed} provided $X$ is closed in any
space $Y$ of the class $C$ containing $X$ as a subspace. It is
well known that when $C$ is the class of Tychonoff spaces than
$C$-closedness coincides with compactness. For the class of
Hausdorff spaces the following conditions for a space $X$ are
equivalent [1, 3.12.5]

(1) The space $X$ is H-closed.

(2) If $\Cal V$ is a centered family of open subsets of $X$ then
$\bigcap \{\ol V:V\in\Cal V\}\not=\0$.

(3) Every ultrafilter in the family of all open subsets of $X$ is
convergent.

(4) Every cover $\Cal U$ of the space $X$ contains a finite
subfamily $\Cal V$ such that $\bigcup\{\ol V:V\in\Cal V\}=X$.

The group $G$ with topology $\tau$ is called a {\it
paratopological group} if the multiplication on the group $G$ is
continuous. If the inversion on the group $G$ is continuous then
$(G,\tau)$ is a {\it topological group}. A group $(G,\tau)$ is
paratopological if and only if the following conditions (known as Pontrjagin
conditions) are satisfied for base $\Cal B$ at unit $e$ of $G$ [4,5].

1. $\bigcap\{UU^{-1}:U\in\Cal B\}=\{e\}$.

2. $(\forall U,V\in\Cal B)(\exists W\in\Cal B):W\subset U\cap V$.

3. $(\forall U\in\Cal B)(\exists V\in\Cal B):V^2\subset U$.

4. $(\forall U\in\Cal B)(\forall u\in U)(\exists V\in\Cal
B):uV\subset U$.

5. $(\forall U\in\Cal B)(\forall g\in G)(\exists V\in\Cal
B):g^{-1}Vg\subset U$.

The paratopological group $G$ is a topological group if and only if

6. $(\forall U\in\Cal B)(\exists V\in\Cal B):V^{-1}\subset U$.

A topological group is {\it absolutely  closed\/} if it is closed
in every Hausdorff topological group containing it as a
topological subgroup. A topological group $G$ is H-closed if and only if it
is {\it Rajkov-complete\/}, that is complete with respect to the
upper uniformity which is defined as the least upper bound $\Cal
L\vee\Cal R$ of the left and the right uniformities on $G$. Recall
that the sets $\{(x,y): x^{-1}y\in U\}$, where  $U$  runs  over  a
base  at  unit  of $G$, constitute a base of entourages for the
left uniformity $\Cal L$ on $G$.  In  the  case  of  the  right
uniformity  $\Cal R$,  the condition $x^{-1}y\in U$ is replaced by
$yx^{-1}\in U$. The {\it Rajkov completion $\hat G$\/} of a
topological group $G$ is the completion of $G$ with respect to the
upper uniformity $\Cal L\vee\Cal R$. For every topological group
$G$ the space $\hat G$ has a natural structure of a topological
group. The group $\hat G$ can be defined as a unique (up to an
isomorphism) Rajkov complete group containing $G$ as a dense
subgroup.

A paratopological group is {\it H-closed} if it is closed in
every Hausdorff paratopological group containing it as a subgroup.
In the present section we shall consider H-closed paratopological
groups.

{\sl Question} Let $G$ be a regular paratopological group
which is closed in every regular paratopological group containing
it as a subgroup. Is $G$ H-closed?

\SecondPage{ON H-CLOSED PARATOPOLOGICAL GROUPS}

\proclaim
{1. Lemma} Let $(G,\tau)$ be a paratopological group.
If there exists a paratopology $\sigma$ on the group $G\times\Z$
such that $\sigma|G\subset\tau$ and $e\in\ol{(G,1)}^\sigma$ then
$(G,\tau)$ is not H-closed.
\endproclaim
{\sl Proof.} We shall build the paratopology $\rho$ on the group $G\times\Z$
such that $\rho|G=\tau$ and $\ol G^\rho\not=G$. Determine the base of unit
$\Cal B_\rho$ as follows. Let $S=\{(x,n):x\in G,n>0\}$. For every neighborhoods
$U_1\in\tau$, $U_2\in\sigma$ such that $U_1\subset U_2$ put $(U_1,U_2)=U_1\cup
(U_2\cap S)$. Put $\Cal B_\rho=\{(U_1,U_2):U_1\in\Cal B_\tau,U_2\in\Cal
B_\sigma\}$. Verify that $\Cal B_\rho$ satisfies the Pontrjagin conditions.

1. It is satisfied since $(U_1,U_2)\subset U_2$.

2. It is satisfied since $(U_1\cap V_1,U_2\cap V_2)
\subset (U_1,U_2)\cap(V_1,V_2)$.

3. Select $V_2\in\Cal B_\sigma$ and $V_1\in\Cal B_\tau$ such that
$V_2^2\subset U_2$, $V_1^2\subset U_1$ and $V_1\subset V_2$. Let
$y_1,y_2\in(V_1,V_2)$. The following cases are possible

A. $y_1,y_2\in V_1$. Then $y_1y_2\in V_1^2\in (U_1,U_2)$.

B. $y_1\in V_1,y_2\in V_2\cap S$. Then
 $y_1y_2\in V_2^2\in U_2$. Since $y_1\in G$ and $y_2\in S$
then $y_1y_2\in S$ and hence $y_1y_2\in U_2\cap S$.

C. $y_1\in V_2\cap S,y_2\in V_1$ is similar to the case B.

D. $y_1,y_2\in V_2\cap S$. Since $S$ is a semigroup then
$y_1y_2\in U_2\cap S$.

4. Let $y\in (U_1,U_2)$. There exist $V_2\in\Cal B_\sigma$
and $V_1\in\Cal B_\tau$ such that $yV_2\subset U_2$ and
$V_1\subset V_2$. The following cases are possible

A. $y\in U_1$. We may suppose that  $yV_1\subset U_1$.
Since $y\in G$ then $y(V_2\cap S)\subset U_2\cap S$.

B. $y\in U_2\cap S$.
Since $V_1\subset G$ then $yV_1\in U_2\cap S$.
Since $S$ is a semigroup and $y\in S$
then $y(V_2\cap S)\subset U_2\cap S$.
Therefore $y(V_1,V_2)\subset (U_1,U_2)$.

5. Let $(g,n)\in G\times\Z$. There exist $V_2\in\Cal B_\sigma$ and
$V_1\in\Cal B_\tau$ such that $V_1\subset V_2$, $g^{-1}V_1g\subset
U_1$ and $g^{-1}V_2g\subset U_2$. Then
$(g,n)^{-1}(V_1,V_2)(g,n)=g^{-1}(V_1,V_2)g=g^{-1}(V_1\cup (V_2\cap
S)g \subset U_1\cup (U_2\cap S) =(U_1,U_2)$.

Therefore $(H,\rho)$ is a paratopological group. Since $(U_1,U_2)\cap G=U_1$
then $\rho|G=\tau$.

Since ${e}\in\ol{(G,1)}^\sigma$ then for every $U_2\in\Cal
B_\sigma$ there exists $g\in G$ such that $(g,1)\in U_2$. Then
$g\in (e,-1)(U_2 \cap S)$ and therefore $(e,-1)\in\ol G^\rho$.\qed

A group topology $\tau_1$ on the group $G$ is called complementable if there
exist a nondiscrete group topology $\tau_2$ on $G$ and neighborhoods
$U_i\in\tau_i$ such that $U_1\cap U_2=\{e\}$. In this case we say that $\tau_2$
is a {\it complement} to $\tau_1$. Proposition 1.4 from [1] implies that in
this case a topology $\tau_1\wedge \tau_2$ is Hausdorff.

A Banach measure is a real function $\mu$ defined on the family of
all subsets of a group $G$ which satisfies the following
conditions:

(a) $\mu(G)=1$.

(b) if $A,B\subset G$ and $A\cap B=\0$ then $\mu(A\cup
B)=\mu(A)+\mu(B)$.

(c) $\mu(gA)=\mu(A)$ for every element $g\in G$ and for every
subset $A\subset G$.

\proclaim {2. Lemma} {\rm [3, p.37]}. Let $G$ be an abelian group and let $\mu$
be a Banach measure on $G$. Let $\tau$ be a group topology on $G$. Suppose that
the set $nG$ is $U$-unbounded for some natural number $n$ and for some
neighborhood $U$ of zero in $(G,\tau)$. Then $\mu(\{x\in G:nx\in gW\})=0$ for
every element $g\in G$ and for every neighborhood $W$ of zero satisfying
$WW^{-1}\subset U$.
\endproclaim

Let $U$ be a neighborhood of zero in a topological group $(G,\tau)$. We say
that a subset $A\subset G$ is $U$-unbounded if $A\not\subset KU$ for every
finite subset $K\subset G$.

Given any elements $a_0,a_1,\dots,a_n$ of an abelian group $G$ put

$$Y(a_0,a_1,\dots,a_n)=\{a_0^{x_0}a_1^{x_1}\cdots a_n^{x_n}:
0\le x_i\le i+1,i\le n,\sum x_i^2>0\},$$

$$X(a_0,a_1,\dots,a_n)=\{a_0^{x_0}a_1^{x_1}\cdots a_n^{x_n}:
-(i+1)\le x_i\le i+1,i\le n\}.$$
Then $X(a_0,a_1,\dots,a_n)=Y(a_0,a_1,\dots,a_n)Y(a_0,a_1,\dots,a_n)^{-1}$.

\proclaim {3. Lemma} Let $(G,\tau)$ be an abelian paratopological group of the
infinite exponent. If there exists a neighborhood $U\in\Cal B_\tau$ such that a
group $nG$ is $UU^{-1}$-unbounded for every natural number $n$ then the 
paratopological group
$(G,\tau)$ is not H-closed.
\endproclaim

{\sl Proof.} Define a seminorm $|\cdot|$ on the group $G$ such that for all
$x,y\in G$ holds $|xy|\le |x|+|y|$. Suppose that there exists a non periodic
element $x_0\in G$. Determine a map $\phi_0:\langle x_0\rangle\to\Z$ putting
$\phi_0(x_0^n)=n$. Since $\Q$ is a divisible group then the map $\phi_0$ can be
extended to a homomorphism $\phi:G\to\Q$. Put $|x|=|\phi(x)|$ for every element
$x\in G$. If $G$ is periodic then put $|e|=0$ and $|x|=[\ln ord(x)]+1$, where
$ord(x)$ denotes the order of the element $x$.

Fix a neighborhood $V\in\Cal B_\tau$ such that 
$V^2\subset U$
 and put
$W=VV^{-1}$. We shall construct a sequence $\{a_n\}$ such that

\item{(a)} $|a_n|>n$.
\item{(b)} $W\cap X(a_0,a_1,\dots,a_n)=\{e\}$.
\item{(c)} $Y(a_0,a_1,\dots,a_n)\not\ni e$.
\item{(d)} if $-n\le k\le n,k\not=0$ then $a_n^k\not\in 2X(a_0,a_1,\dots,a_{n-1})$.

Take any element $a_0\not\in W$. Suppose that the elements
$a_0,\dots,a_n$ have been chosen satisfying conditions (a) and (b).
Put $$B_n=\{x\in G:(\forall g\in X(a_0,a_1,\dots,a_{n-1}))
(\forall k\in\Z\bs\{0\} :-e^{n+1}\le k\le e^{n+1}):kx\not\in gW\}.$$

If the group $ G$ is periodic then $|x|>n$ for every element $x\in B_n$. Lemma
2 implies that $\mu(B_n)=1$. If the group $ G$ is not periodic then the
construction of the seminorm $|\cdot|$ implies that $\mu(\{x\in G:|x|\le
n\})=\mu(\phi^{-1}[-n;n])=0$. In both cases there exists an element $a_n\in
B_n$  such that $|a_n|>n$. Then $W\cap X(a_0,a_1,\dots,a_n)=\0$. Considering a
subsequence and applying condition (a) we can satisfy conditions (c) and (d)
also.

Define a base $\Cal B_{\tau\{a_n\}}$ at the unit of group topology
$\tau\{a_n\}$ on the group $ G\times\Z$ as follows. Put
$A_n^+=\{(e,0)\}\cap\{(a_k,1):k\ge n\}$. For every increasing
sequence $\{n_k\}$ put $A[n_k]=\bigcup\limits_{l\in\N}
A_{n_1}^+\cdots A_{n_l}^+$. Put $\Cal B_{\tau\{a_n\}}=\{A[n_k]\}$.
We claim that $(G\times Z,\tau\{a_n\})$ is a zero dimensional
paratopological group.

Put $F=\bigcup\limits_{n\in\omega}X(a_0,a_1,\dots,a_n)$. Let
$A[n_k]\in \Cal B_{\tau\{a_n\}}, (x,n_x)\not\in A[n_k]$. If
$x\not\in F$ then $(x,n_x)A[n]\cap A[n_k]=\0$. Let $x\in
X(a_0,a_1,\dots,a_m)$. Put $m_k=m+k$. Suppose that
$(x,n_x)A[m_k]\cap A[n_k]\not=\0$. Select the minimal $k$ such
that $(x,n_x)(A_{m_1}^+\cdots A_{m_k}^+)\cap A[n_k]\not=\0$. Let
$$(x,n_x)(a_{l_1},1)\cdots (a_{l_k},1)=(a_{l'_1},1)\cdots
(a_{l'_{k'}},1)\leqno{(*)}$$ and for all $i,i'$ holds $m_i\le
l_i\le l_{i+1}$, $n_i'\le l'_{i'}\le l'_{i'+1}$. Remark that a
member $a_q$ occurs in each part of the equality $(**)$ no more
than $q$ times. If $l_k>l'_{k'}$ then if we move all members which
are not equal to $(a_{l_k},1)$ from the left side of the equality
(*) to the right one, we obtain contradiction to condition $(d)$.
The case $l_k<l'_{k'}$ is considered similarly. Therefore
$l_k=l'_{k'}$, a contradiction to that $k$ is the minimal number
such that the equality (*) holds. It is showed similarly that if
$x\not=e$ and $m_k=m+k+1$ then $(x,n_x)\not\in A[m_k]$. If $x=e$
and $n_x\not=0$ then condition (c) implies that $A[n]\not\ni
(x,n_x)$. Hence Pontrjagin condition 1 for $\Cal B_{\tau\{a_n\}}$
is satisfied. Since $A[n_{2k}]^2\subset A[n_k]$, Pontrjagin
condition 3 is satisfied. All other Pontrjagin conditions are
obvious.

Condition (b) implies that $A[n]A[n]^{-1}\cap VV^{-1}=\{(e,0)\}$. Therefore the
topology $\tau{\{a_n\}}_g$ is a complement to the topology
$(\tau\times\{0\})_g$, where $\tau\times\{0\}$ is the product topology on the
group $(G,\tau)\times\Z$. Therefore the topology
$\sigma=\tau{\{a_n\}}(\tau\times\{0\})$ is Hausdorff. Since
$(e,0)\in\ol{(G,1)}^{\tau{\{a_n\}}}\subset\ol{(G,1)}^\sigma$ then $(G,\tau)$ is
not H-closed.\qed

We shall need the following lemma.

\proclaim{4. Lemma} Let $G$ be a paratopological group and $H$ be a normal
subgroup of the group $G$. If $H$ and $G/H$ are topological groups then $G$ is
a topological group.
\endproclaim
{\sl Proof.} Let $U$ be an arbitrary neighborhood of the unit. There exist
neighborhoods $V,W$ of the unit such that $V\subset U$, $(V^{-1})^2\cap
H\subset U$ and $W\subset V$, $W^{-1}\subset VH$. If $x\in W^{-1}$ then there
exist elements $v\in V,h\in H$ such that $x=vh$. Then $h=v^{-1}x\in
V^{-1}W^{-1}\cap H\subset U$. Therefore $x\in VU\subset U^2$. Hence $G$ is a
topological group.\qed

\proclaim{5. Theorem} An abelian topological group $(G,\tau)$ is H-closed if
and only if
$(G,\tau)$ is Rajkov complete and for every group topology $\sigma\subset\tau$
on $G$ the quotient group $\hat G/G$ is periodic, where $\hat G$ is the Rajkov
completion  of the group $(G,\sigma)$.
\endproclaim
{\sl Proof.} Suppose that there exists a group topology $\sigma\subset\tau$ on
$G$ such that the quotient group $\hat G/G$ is not periodic, where $\hat G$ is
the Rajkov completion of the group $(G,\sigma)$. Select a non periodic element
$x\in\hat G$ such that $\langle x\rangle\cap G=\{e\}$. Then $G\times \langle
x\rangle$ is naturally isomorphic to a group $G\times\Z$ and Lemma 1 implies
that the group $(G,\tau)$ is not H-closed.

Let a paratopological group $(H,\tau')$ contains $(G,\tau)$ as non closed
subgroup. Since $G$ is abelian then $\ol G$ is an abelian semigroup. Choose an
arbitrary element $x\in\ol G\bs G$. Then a group hull $F=\langle G,x\rangle$
with a topology $\tau'|F$ is an abelian paratopological group. Then the group
$G$ is dense in $(F,{\tau'}_g)$. Since the Rajkov completion $\hat F$ of the
topological group $(F,\tau'|F_g)$ is periodic then there exists a natural
number $n$ such that $x^n\in G$. Therefore $F^n\subset G$. Lemma 4 implies that
$F$ is a topological group and therefore $G$ is closed in $(F,{\tau'}_g)$, a
contradiction.\qed

\proclaim{6. Corollary} A Rajkov completion of a isomorphic condensation of
H-closed abelian topological group is H-closed.
\endproclaim

\proclaim{7. Proposition} Let $G$ be a Rajkov complete topological group, $H$
be H-closed paratopological subgroup of the group $G$. If a group $G/H$ has
finite exponent then $G$ is an H-closed paratopological group.
\endproclaim
{\sl Proof.} Select a number $n$ such that $g^n\in H$ for every element $g\in
G$. Let $F\supset G$ be a paratopological group. Since $H$ is closed in $F$
then for every element $g\in\ol G$ holds $g^n\in H$. Denote continuous maps
$\phi:\ol G\to\ol G$ as $\phi(g)=g^{n-1}$ and $\psi:\ol G\to H$ as
$\psi(g)=(g^{n})^{-1}$. Then for every element $g\in\ol G$ holds
$g^{-1}=\phi(g)\psi(g)$ and therefore the inversion on the group $\ol G$ is
continuous. Since $\ol G$ is a topological group and $G$ is Rajkov complete
then $\ol G=G$.\qed

\proclaim{8. Proposition} Let $G$ be a paratopological group and $K$ be a
compact normal subgroup of the group $G$. If a group $G/K$ is H-closed then the
group $G$ is H-closed.
\endproclaim
{\sl Proof.} Suppose that there exists a paratopological group $F$ containing
the group $G$ such that $\ol G\not=G$. Since $K$ is compact then $F/K$ is a
Hausdorff paratopological group by Proposition 1.13 from [4]. Let $\pi:F\to
F/K$ be the standard map. Then
$\ol{G/K}\supset\pi(\ol{\pi^{-1}(G/K)})\supset\pi(\ol G)\not=\pi(G)=G/K$. This
implies that the group $G/K$ is not H-closed, a contradiction.\qed

Let $G$ be a topological group, $N$ be a closed normal subgroup of the group
$G$. Then if $N$ and $G/N$ are Rajkov complete so is the group $G$ [5]. This
suggests the following

{\sl 9. Question.} Let $G$ be a paratopological group, $N$ be a closed normal
subgroup of the group $G$ and the groups $N$ and $G/N$ are H-closed. Is the
group $G$ H-closed?

Let $(G,\tau)$ be a paratopological group. Then there exists the finest group
topology $\tau_g$ coarser than $\tau$ (see [2]), which is called {\it the group
reflection} of the topology $\tau$.

\proclaim{10. Proposition} Let $(G,\tau)$ be an abelian paratopological group.
If $(G,\tau_g)$ is H-closed then $(G,\tau)$ is H-closed. If $(G,\tau)$ is
H-closed and $(G,\tau_g)$ is Rajkov complete then $(G,\tau_g)$ is H-closed.
\endproclaim
{\sl Proof.} Suppose that the group $(G,\tau_g)$ is H-closed and $(G,\tau)$ is
not. Let a paratopological $(H,\hat\tau)$ contains $(G,\tau)$ as non closed
subgroup. Without loss of generality we may suppose that there exists an
element $x\in H\bs G$ such that $H=\langle G,x\rangle$ and the group $H$ is
abelian. Let $\hat\tau_g$ be the group reflection of the topology $\hat\tau$.
Since $\hat\tau_g|G\subset\tau_g$ then Theorem 5 implies that the group $H/G$
is periodic. Without loss of generality we may suppose that $x^p\in G$ for some
prime $p$.

Denote the family of neighborhoods at unit in the topology $\tau$
as $\Cal B_{\hat\tau}$. Let $U\in\Cal B_{\hat\tau}$. If $U\cap
xG=\0$ then there exists a neighborhood $V$ of unit such that
$V^p\subset U$ and thus $V\subset G$ and $G$ is open in
$(H,\hat\tau)$. Therefore a set $\Cal F=\{x^{-1}(xG\cap
U):U\in\Cal B_\tau\}$ is a filter. Let $U\in\Cal B_{\hat\tau}$.
There exists $V\in\Cal B_{\hat\tau}$ such that $V^p\subset U$.
Then $(xG\cap V)^p\subset U$. Let $xg\in (xG\cap V)$. Then
$x^{-1}(xG\cap V)\subset x^{-1}((xg)^{1-p}(xG\cap V)^p)\cap G
\subset x^{-p}g^{1-p}(U\cap G)$ and hence $\Cal F$ is a Cauchy
filter in the group $(G,\tau_g)$. Let $h\in G$ be a limit of the
filter $\Cal F$ on the group $(G,\tau_g)$. But then for every
neighborhood of the unit $U$ in the topology $\hat\tau_g$ holds
$U\cap xhU\supset U\cap xh(U\cap G)\not=\0$ and therefore
$(H,\hat\tau_g)$ is not Hausdorff, a contradiction.

Let $(G,\tau_g)$ is Rajkov complete and $(G,\tau_g)$ is not H-closed. Then
Theorem 5 implies that there exists a group topology $\sigma\subset\tau$ on $G$
such that the quotient group $\hat G/G$ of the Rajkov completion $\hat G$ of
the group $(G,\sigma)$ is not periodic. Then Lemma 1 implies that a group
$(G,\tau)$ is not H-closed.\qed

\proclaim{11. Lemma} Let topological group $(H,\sigma_H)$ be a closed subgroup
of an abelian topological group $(G,\tau)$ and $\sigma_H\subset\tau|H$. Then
there exists a group topology $\sigma\subset\tau$ on the group $G$ such that
$\sigma|H=\sigma_H$.
\endproclaim
{\sl Proof} Let $\Cal B_\tau$ and $\Cal B_{\sigma H}$ be bases of unit of
$(G,\tau)$ and $(H,\sigma_H)$ respectively.

Put $\Cal B_\sigma=\{U_1U_2:U_1\in\Cal B_\tau,U_2\in\Cal B_{\sigma
H}\}$. Verify that the family $\Cal B_\sigma$ satisfies the
Pontrjagin conditions.

2. It is satisfied since $(U_1\cap V_1)(U_2\cap V_2)\subset
U_1U_2\cap V_1V_2$.

3. Select $V_2\in\Cal B_{\sigma H}$ and $V_1\in\Cal B_\tau$ such
that $V_2^2\subset U_2$, $V_1^2\subset U_1$. Then
$(V_1V_2)^2\subset U_1U_2$.

4. Let $y\in U_1U_2$. Then there exist points $y_1\in U_1$ and
$y_2\in U_2$ such that $y=y_1y_2$. Therefore there exist a
neighborhoods $V_1\in\Cal B_\tau$ and $V_2\in\Cal B_{\sigma H}$
such that $y_iV_i\subset U_i$. Then $yV_1V_2\subset U_1U_2$.

5. It is satisfied since $G$ is abelian.

6. $(U_1^{-1}U_2^{-1})^{-1}\subset U_1U_2$.

1. Since all others Pontrjagin conditions are satisfied then it
suffice to show that $\bigcap\Cal B_\sigma=\{e\}$. Let $x\in G$
and $x\not=e$. If $x\in H$ then there exists $U_2\in\Cal B_{\sigma
H}$ such that $U_2^2\not\ni x$ and $U_1\in\Cal B_\sigma$ such that
$U_1\cap H\subset U_2$. Then $U_1U_2\cap\{x\}=U_1U_2\cap\{x\}\cap
H\subset U_2^2\cap\{x\}=\0$. If $x\not\in H$ then $(G\bs
xH)H\not\ni x$.

Therefore $(G,\sigma)$ is a topological group. Since $U_1U_2\cap H=(U_1\cap
H)U_2$ then $\sigma|H=\sigma_H$. \qed

\proclaim{12. Proposition} A closed subgroup of an H-closed abelian group is
H-closed.
\endproclaim
{\sl Proof.} Let $H$ be a closed subgroup of an H-closed abelian group
$(G,\tau)$. Then $G$ and $H$ are Rajkov complete. Let $\sigma_H\subset\tau|H$
be a group topology on the group $H$. Lemma 11 implies that there exists a
group topology $\sigma$ on the group $G$ such that $\sigma|H=\sigma_H$. Let
$(\hat G,\hat\sigma)$ be the Rajkov completion of the group $(G,\sigma)$. Then
a closure $\ol H^{\hat\sigma}$ of the group $H$ in the group $(\hat
G,\hat\sigma)$ is a Rajkov completion of the group $(H,\sigma_H)$. Let
$x\in\ol{H}^{\hat\sigma}$. Theorem 5 implies that there exists $n>0$ such that
$x^n\in G$. Since $\ol{H}^{\hat\sigma}\cap G=H$ then $x^n\in H$. Therefore
Theorem 5 implies that $H$ is H-closed.\qed

\proclaim{13. Proposition} Let $G$ be a H-closed abelian topological group.
Then $K=\bigcap_{n\in\N}\ol{nG}$ is compact and for each neighborhood $U$ of
zero in $G$ there exists a natural $n$ with $\ol{nG} \subset KU$.
\endproclaim
{\sl Proof.} Let $\Phi$ be a filter on $G$ with a base $\{\ol{nG}:n\in\N\}$,
and $\Psi$ be an arbitrary ultrafilter on $G$ with $\Psi\supset\Phi$. Let $U$
be a closed neighborhood of the unit in $G$. Lemma 2 implies that there exists
a number $n$ such that the set $\ol{nG}$ is $U$-bounded. Since $\ol{nG}\in\Phi$
and $\Psi$ is an ultrafilter, there exists  $g\in G$ with $gU\in\Psi$. Hence
$\Psi$ is a Cauchy filter on $G$. By the completeness of $G$, $\Psi$ is
convergent. Therefore each ultrafilter $\Psi$ on $G$ with $\Psi\supset\Phi$
converges. In particular each ultrafilter on $K$ is convergent, and since $K$
is closed, $K$ is compact.

To show that there exists a number $n$ with  $\ol{nG} \subset KU$ it suffices
to prove that $KU\in\Phi$. Assume that $KU\not\in\Phi$. Then there exists an
ultrafilter $\Psi\supset\Phi$ with $G\bs KU\in\Psi$. As we have proved, $\Psi$
is convergent. Clearly $\lim\Psi\in K$. Therefore $KU\in\Psi$ which is a
contradiction. Hence $KU\in\Phi$, and this completes the proof. \qed

\proclaim{14. Corollary} A divisible abelian H-closed topological group is
compact.\qed
\endproclaim

\proclaim{15. Proposition} Every H-closed abelian topological group is a union
of compact groups.
\endproclaim
{\sl Proof.} Let $G$ be such a group. It suffice to show that every element
$x\in G$ is contained in a compact subgroup. Let $X$ be the smallest closed
subgroup of $G$ containing the element $x$. Then $X=\bigcup_{k=0}^n
(kx+\ol{nX})$ for every natural $n$. Let $U$ be an arbitrary neighborhood of
the zero. By Lemma 15 there exists a natural number $n$ such that $nG$ is
$U$-bounded. Then $X$ is also $U$-bounded. Hence $X$ is a precompact group.
Since $X$ is Rajkov complete then $X$ is compact.\qed

{\bf 16. Conjecture.} An abelian topological group $G$ is H-closed if 
and only if $G$ is
Rajkov complete and $nG$ is precompact for some natural $n$.

\proclaim{17. Proposition} The Conjecture 16 is true provided the group
$(G,\tau)$ satisfies the following two conditions:

(1) There exists a $\sigma$-compact subgroup $L$ of $G$ such that
$G/L$ is periodic.

(2) There exists a group topology $\tau'\subset\tau$ such that
the Rajkov completion $\hat G$ of the group $(G,\tau')$ is Baire.
\endproclaim
{\sl Proof.} Let $G$ be such a group and $L=\bigcup_{k\in\N} L_k$ be a union of
compact subsets $L_k$. Put $G(n,k)=\{x\in\hat G:nx\in L_k\}$ for every natural
$n$ and $k$. Then every set $G(n,k)$ is closed. By Theorem 5 holds $\hat
G=\bigcup_{n,k\in\N} G(n,k)$. Since $\hat G$ is Baire then there exist natural
numbers $n$ and $k$ such that $\inte G(n,k)\not=\0$. Then $F=G(n,k)-G(n,k)$ is
a neighborhood of the zero.  By Corollary 6 the group $\hat G$ is H-closed. Put
$K=\bigcap_{n\in\N}\ol{n\hat G}$. By Proposition 13 there exists a natural $m$
such that $m\hat G\subset F+K$. Then $mnG\subset mn\hat G\subset L_k-L_k+K$ and
hence the group $mnG$ is precompact.\qed

\centerline{------------------------------} 

\vskip15pt
\parskip3pt

\item{1.} {\it Engelking R.} General topology. --  Monografie
Matematyczne, Vol. 60, Polish Scientific Publ. -- Warsaw, 1977.

\item{2.} {\it Graev M.I.} Theory of topological groups //
UMN, 1950 (in Russian).

\item{3.} {\it Protasov I., Zelenyuk E.}, Topologies on Groups
Determined by Sequences. --  VNTL Publishers, 1999.

\item{4.} {\it Ravsky O.V.} // Paratopological groups I //
Matematychni Studii. -- 2001. -- Vol.~16. -- \char 194 1. -- P.37-48.

\item{5.} {\it Ravsky O.V.} // Paratopological groups II //
Matematychni Studii. -- 2002. -- Vol.~17. -- \char 194 1. -- P.93-101.
\bye
\vskip15pt

\centerline{\bf ПРО H-ЗАМКНЕНI ПАРАТОПОЛОГIЧНI ГРУПИ}
\vskip2pt
\centerline{\bf О. Равський}
\vskip4pt
\centerline{\it Львўвський нацўональний унўверситет ўменў Ўвана Франка}
\centerline{\it вул.Унўверситетська, 1 79000 Львўв, Укра•на}
\vskip15pt 

Гаусдорфова паратопологiчна група називаїться H-замкненою, якщо вона замкнена у
довiльнiй гаусдорфовiй паратопологiчнiй групi, що •• мiстить. Отримано критерiй
H-замкненостi абелево• топологично• групи i для деяких класiв абелевих
паратопологiчних груп отримано простi критерi• H-замкненостi.

\noindent
Ключовў слова: {\it паратопологўчна група, мўнўмальна топологўчна 
група, абсолютно зам\-кне\-на то\-по\-ло\-гўч\-на група.}




\bye